\documentclass[reqno]{amsart}

\usepackage{graphicx}
\usepackage{color}
\newtheorem{theorem}{Theorem}
 \newtheorem{lem}{Lemma}
  \newtheorem{remark}{Remark}

\UseRawInputEncoding
\usepackage[russian, english]{babel}
 \usepackage{setspace}
\def\dfrac#1#2{\displaystyle{#1\over #2}}

\def\bv{{\bf V}}
\def\bV{{\bf V}}

\def\Div{\mbox{div}\,}

\def\bx{{\bf x}}
\def\bE{{\bf E}}
\renewcommand{\thefootnote}{\fnsymbol{footnote}}

\begin{document}

\title[
 Oscillations of a Collisional  Plasma]{
On Radially Symmetric Oscillations of a Collisional  Cold Plasma
}


\author{Olga S. Rozanova* \, and \,  Maria I. Delova}




\renewcommand{\thefootnote}{\fnsymbol{footnote}}
\footnotetext{\emph{* \,Corresponding author:} Olga Rozanova}
\renewcommand{\thefootnote}{\arabic{footnote}}

\maketitle
\begin{center}
{\it Lomonosov Moscow State University, \\ MSU, Main building, Moscow, 119991, Russian Federation, \\e-mail: rozanova@mech.math.msu.su }
\end{center}

\begin{abstract}
We study the influence of the friction term on the radially symmetric solutions of the repulsive
Euler-Poisson equations with a non-zero background, corresponding to
cold plasma oscillations in many spatial dimensions. It is shown that for any arbitrarily small  non-negative constant friction coefficient, there exists a neighborhood of the zero equilibrium in the $C^1$ norm such that the solution of the Cauchy problem with initial data belonging to this neighborhood remains globally smooth in time. Moreover, this solution stabilizes to  zero as $t\to\infty$. This result contrasts with the situation of zero friction, where any small deviation from the zero equilibrium generally leads to a blow-up. Our method  allows us to estimate the lifetime of smooth solutions. Further, we prove that for any initial data, one can find such  coefficient of friction that the respective solution to the Cauchy problem keeps smoothness for all $t>0$ and stabilizes to zero.
We also present the results of numerical experiments for  physically reasonable situations, which allows us to estimate the value of the friction coefficient, which makes it possible to suppress the formation of singularities of solutions.
\end{abstract}

\section{Introduction}

We study a frictional version of the repulsive Euler-Poisson
equations
\begin{eqnarray}\label{EP}
\dfrac{\partial n }{\partial t} + \Div(n \bv)=0,\quad
\dfrac{\partial \bv }{\partial t} + \left( \bv \cdot \nabla \right)
\bv =\,k \,  \nabla \Phi\, -\,\nu \,\bv, \quad \Delta \Phi =n-n_0,
\end{eqnarray}
where the  scalar functions $n$ and $\Phi$ are the density and a repulsive (for $k>0$) force potential, respectively,   the vector $\bv$ is the velocity, they depend on
the time $t$ and the point $x\in {\mathbb R}^d $, $d\ge 1$. Here
$n_0={\rm const}> 0$ is the density background, $\nu={\rm const}>0$ is the friction coefficient.

If we denote $\nabla \Phi = -\bE$, and set  $n_0=1$, such that
\begin{eqnarray}
n=1- \Div \bE,\label{n}\end{eqnarray} we can remove $ n $
from \eqref{EP} and
 rewrite it as
\begin{eqnarray}\label{4}
\dfrac{\partial \bv }{\partial t} + \left( \bv \cdot \nabla \right)
\bv = \, - \bE - \nu \bv,\quad \frac{\partial \bE }{\partial t} + \bv \Div \bE
 = \bV.
\end{eqnarray}

In this paper we study the Cauchy problem for \eqref{EP} or  \eqref{4} and our main concern is to study initial data that guarantee a globally smooth
solution.

System \eqref{4} corresponds to the  hydrodynamics of ``cold" or electron plasma in the
non-relativistic approximation in dimensionless quantities  (see, e.g.,~\cite{ABR78}, \cite{david72},  \cite {GR75}).
In this interpretation the friction coefficient characterizes the intensity of electron-ion collisions during plasma oscillations.
The cold plasma equations is now very popular object of study, may be more popular than the Euler-Poisson equations themselves.
The reason is that the cold plasma in used in the accelerators of electrons in the
wake wave of a powerful laser pulse \cite {esarey09}.  From this point of view the initial data (i.e. the initial laser pulse) that corresponds to a solution that cannot survive being smooth are not applicable technically.

System \eqref{4} in 1D case was studied \cite{RChD20}, however, in the multidimensional case it is very difficult from both mathematical and physical points of view. Indeed,  it describes a non-hyperbolic superposition of different types of waves, each of them have a tendency to break out in a finite time. Therefore the theoretical results here are very scarce.

The situation is more optimistic if we restrict ourselves to the class of axisymmetric solutions. Thus, we consider one-dimensional solutions, given on the half-line. In \cite{ELT}, \cite{CT}, \cite{WTB}, \cite{T21} it was shown that for $n_0=0$, $k>0$ and $n_0\ge 0$, $k<0$ in the non-frictional case there is a threshold in terms of the initial data. Namely, one can specify exactly the class of initial data corresponding to a globally smooth solution, and these data form a neighborhood of the zero equilibrium in the $C^1$ -norm. As it has been recently shown \cite{R22_Rad}, for $n_0>0$, $k>0$, $\nu=0$ the situation is strikingly different: namely, for $ d\ne 1$ and $d\ne 4$ an arbitrarily small perturbation of the zero equilibrium blows up in the general case. The exception is the initial data in the form of a simple wave, starting from which the solution can remain globally smooth and tend to an affine solution as $t\to \infty$. In any case, the initial data corresponding to simple waves form a zero-measure manifold in the neighborhood of the zero equilibrium.

In this paper, we study the effect of constant coefficient of friction on the blow-up process. Namely, we establish that the presence of friction "normalizes" the situation with the threshold for the initial data. Namely, for an arbitrarily small $\nu>0$ and any $d$, there exists a neighborhood of the zero stationary state in the $C^1$-norm such that the corresponding solution of the Cauchy problem preserves smoothness (Theorem 1). For small $\nu$, Theorem 2 gives sufficient conditions guaranteeing blow-up or  non-blow-up in terms of initial data, which can be applied to numerical tests.
 Besides, we show that for any initial data, one can find $\nu$ such that the corresponding solution of the Cauchy problem is globally smooth (Theorem 3). In other words, this situation is absolutely analogous to $d=1$, and the increase in spatial dimension does not lead to any new phenomena.

Thus, we consider  radially symmetric solutions of  \eqref{4}
\begin{eqnarray}\label{sol_form}
\bV=F(t,r)x, \quad \bE=G(t,r)x,
\end{eqnarray}
where $x = (x_1,x_2,...,x_d)$ is the radius-vector, $r=\sqrt{x_1^2+x_2^2+...+x_d^2}$.

The initial data that correspond to these solutions are
\begin{equation}\label{CD1}
(\bv, \bE) |_{t=0}=(\bv_0(r), \bE_0(r))= (F_0(r) {x}, G_0(r) {x} ), \quad (F_0(r)
, G_0(r) ) \in C^2(\bar {\mathbb R}_+).
\end{equation}
We assume that $(\bv_0(r), \bE_0(r))$ are bounded together with their derivatives uniformly on $r\in \bar {\mathbb R}_+$ and denote
 $ \|f\|_{C^1( \overline{ \mathbb{R}}_+)}=\sum\limits_{i=0}^1 \sup\limits_{r \in  \overline{ \mathbb{R}}_+}|f^{(i)}(r)|.$

The
physically natural condition $n|_{t=0}>0$ dictates ${\rm div}\bE <1$, see \eqref{n}.

The main results of the paper are as follows.

\begin{theorem}\label{T1}
For arbitrary small $\nu>0$ there exists  $\varepsilon(\nu)>0$, such that the solution of the problem  \eqref{4} with the initial data \eqref{CD1}, satisfying
\begin{equation}\label{T1CD}
\| \bV_0(r),{\bf E}_0(r)\|_{C^1( \overline{ \mathbb{R}}_+)}<\varepsilon,
\end{equation}
keeps $C^1$ - smoothness  for all $t>0.$ Moreover,
\begin{equation}\label{T1ass}
\| \bV,{\bf E}\|_{C^1( \overline{ \mathbb{R}}_+)}\le {\rm const} \, e^{-\frac{\nu}{2}t} \to 0, \quad t\to \infty.
\end{equation}
\end{theorem}

\bigskip
The next Theorem 2 gives more information about the size of the neighborhood of the origin containing globally smooth solutions in the case of small $\nu$.

To formulate it,
for any fixed $r\in \overline{\mathbb R}_+$  we introduce the following notation.

1. A function that can be determined using the solution $(G(t,r), F(t,r))$ of the auxiliary system \eqref{FG4} under the initial conditions $(G_0(r), F_0 (r))$
\begin{eqnarray}
&&\phi(t,r)=- \frac{d+2}{2} G +(d-2) \nu  F-\frac{(d-2)(d-4)}{2} F^2. \nonumber
\end{eqnarray}

2. Functions of $r$:
\begin{eqnarray}
&&u_0(r)=\Div \bV_0-d F_0, \quad v_0(r)=\Div \bE_0-d G_0, \nonumber\\
&&H_0(r)=u_0,\quad H_1(r)=\left(\frac{d-2}{2}F_0-\frac{\nu}{2}\right)-v_0, \nonumber\\
&& J_+ (r)=1-\frac{\nu^2}{4}- \frac{d+2}{2} G_- +\nu (d-2) F_++(1-\delta) \frac{(d-2)(d-4)}{2} F_+^2,\nonumber\\
&& M_\pm(r)=\left(\frac{1-d G_\mp}{1-d G_0} \right)^\frac{d+2}{2 d}, \quad 0<M_-<M_+,\nonumber
\end{eqnarray}
where 
$G_-=G_-(r)<0$ and $G_+=G_+(r)>0$, $G_+<\frac{1}{d}$ are  the left and right roots of equation  \eqref{Int_2} ($d\ne 2$) or \eqref{Int_d} ($d= 2$), $F=0$ (they depend on $(G_0(r), F_0(r))$), $F_+=F_+(r)$ is given as \eqref{Fp}, $\delta=1$ for $d=3$ and  $\delta=0$, otherwise.

\begin{theorem}\label{T21} Let $\nu<2$.

 a) A sufficient condition on initial data \eqref{CD1} that guaranties the smoothness of the solution of the problem  \eqref{4}, \eqref{CD1}  for all $t>0$ is the following:
\begin{eqnarray}\label{suffcond1}
&&\inf_{r\in  \overline{ \mathbb{R}}_+)} \mathcal F_1 (\nu, \bV_0(r),{\bf E}_0(r))<1,\\\nonumber
&&\mathcal F_1 (\nu, \bV_0,{\bf E}_0)= \frac{2}{\nu} M_+\sqrt{H_0^2+\left(1-\frac{\nu^2}{4}\right)^{-1} H_1^2} \quad e^{\int\limits_0^\infty \|\phi(\tau, .)\|_{L^\infty({\mathbb R}_+)} d\tau}.
\end{eqnarray}

b) If there exists $T>0$ such that
\begin{eqnarray}\label{suffcond2}
&&\inf_{r\in  \overline{ \mathbb{R}}_+)} \mathcal F_2 (T,\nu, \bV_0(r),{\bf E}_0(r))<1,\\
&&\mathcal F_2 (T, \nu, \bV_0,{\bf E}_0)=
 \frac{2}{\nu} M_+\sqrt{H_0^2+\left(1-\frac{\nu^2}{4}\right)^{-1} H_1^2} \quad e^{\left(J_+-1+\frac{\nu^2}{4}\right)T},\label{sc22}
\end{eqnarray}
then
the solution of the problem  \eqref{4}, \eqref{CD1}  preserves smoothness for  $t\in [0,T]$.

c) If the initial data \eqref{CD1} are such that there exists a point $r\in  \overline{ \mathbb{R}}_+$ for which
condition
\begin{eqnarray}\label{suffcond3}
&& \mathcal F_3 (\nu, \bV_0(r),{\bf E}_0(r))\ge 1,\\\nonumber
&&\mathcal F_3 (\nu, \bV_0(r),{\bf E}_0(r))= \frac{2}{\nu} M_- \sqrt{H_0^2+J_+^{-1} H_1^2},
\quad H_0\le 0, \quad H_1<0,\nonumber
\end{eqnarray}
holds. Then the solution  of  problem  \eqref{4}, \eqref{CD1} blows up  within   $t<\frac{ \pi}{\sqrt{J_+}}$.
\end{theorem}

\bigskip
\begin{remark}\label{Rem1}
Condition a) of Theorem 2 guarantees the global existence in $t$ of a smooth solution of the Cauchy problem
      \eqref{4}, \eqref{CD1},  but not explicitly, since it is necessary to solve system \eqref{FG4}, generally speaking, numerically. Condition b) is rougher, but it uses only the initial data. It is easy to see that
      $J_+-1+\frac{\nu^2}{4}>0$ for all possible initial data, so the exponential factor in \eqref{sc22} increases
      with $T$ and the sufficient condition \eqref{suffcond2} cannot be satisfied for all $T>0$.

       However, if we stop trying to obtain explicit sufficient conditions, we can use a numerical procedure to obtain not only sufficient conditions, but even a (numerical) criterion that the solution of the Cauchy problem remains smooth, see Sec.\ref{Example}.
       \end{remark}

\bigskip

 \begin{theorem}\label{T22}
For arbitrary initial data \eqref{CD1} there exists  such $\nu>2$ that  the solution of  problem  \eqref{4}, \eqref{CD1}
keeps $C^1$ - smoothness  for all $t>0$ and the asymptotic property
\begin{equation}\label{T2ass}
\| \bV,{\bf E}\|_{C^1( \overline{ \mathbb{R}}_+)}\le {\rm const} \, e^{-\frac{\nu-\sqrt{\nu^2-4}}{2}t} \to 0, \quad t\to \infty.
\end{equation}
 holds.
\end{theorem}

\bigskip

Theorems \ref{T1}, \ref{T21} and \ref{T22} can be reformulated in the terms of the Euler-Poisson
equations \eqref{EP}. The stationary stationary state in this case is
\begin{equation*}
\bV =0, \quad \Phi = {\rm const},\quad n=1.
\end{equation*}

In this work we use the technique of linearization by means of the Radon lemma, the same as in \cite{R22_Rad}. It turn out to be very convenient for the analysis of the non-strictly hyperbolic systems often arising when studying the reduced cold plasma equations.

The paper is organised as follows: Sec.2 is devoted to auxiliary results on the behavior of solution and its derivatives, Secs.3, 4 and 5 contain the proofs of Theorems 1, 2, and 3, respectively. Sec.6 provides examples of obtaining numerical criteria based on the method of this article for physically reasonable situations.
Sec.7 is devoted to a discussion on the importance of the results for physics and the formulation of future problems in this area.

\section{Behavior of solutions along characteristics}
We use the fact that $\bV=F(t,r)x$, $ {\bf E}=G(t,r)x$ and get
\begin{eqnarray}
\label{FGr}
\frac{\partial F}{\partial t}+Fr \frac{\partial F}{\partial r}=-F^2-G-\nu F, \quad \frac{\partial G}{\partial t}+Fr \frac{\partial G}{\partial r}=F-dFG,
\end{eqnarray}
$$ (\bV^0, {\bf  E}^0)= (F(0,r)x,G(0,r)x)=(F_0(r)x,G_0(r)x), \quad (F_0(r),G_0(r))\in C^2(\bar{\mathbb{R}}_+).$$

\subsection{Physical constraints on solution components}
Let us fix $r_0\in \bar{\mathbb R}_+$. Along the characteristic of  system \eqref{FGr},
 $$\dot{r}=\dfrac{d r}{d t}=Fr, \quad r(0)=r_0,$$ we get
  $$r(t)=r_0 \exp{\int\limits_0^t \, F(\tau) dt}\ge 0, \quad t\ge 0,$$
where $F$ is bounded (see below, Sec.\ref{Bofs}).
Thus, for $r_0 > 0$ we have $r>0$.

Further, along the characteristic the functions $F$ and $ G$ obey the system
\begin{eqnarray}
\label{FG4}
 \dot F=-F^2-G-\nu F, \quad \dot G=F-dFG.
 \end{eqnarray}
Therefore
$$ \frac{\mathrm d G}{(1-d G)\,F}=\frac{\mathrm d r}{Fr},$$
and
\begin{eqnarray}
\label{G}
1-d G=const \cdot r^{-d}.
\end{eqnarray}
  Thus, the sign of the expression $1-d G$ does not change, i.e. $${\rm sign} (1-d G)={\rm sign} (1-d G(0,r_0)).$$ Thus, the motion on the phase plane $(G,F)$, corresponding to  system \eqref{FG4} occurs either in the half-plane $G<\frac{1}{d}$, or in the half-plane $G>\frac{1}{ d}$, or on the line $G=\frac{1}{d}$.

The equilibria of  \eqref{FG4}  are the following:
\begin{itemize}
\item  if $\nu<\frac{2}{\sqrt{d}}$, then there exists the only point
$(F=0, G=0)$, a stable focus;
\item if $\nu=\frac{2}{\sqrt{d}}$, then there exist two points:
$(F=0, G=0)$, a stable focus, and  $(F=-\frac{\nu}{2}, G=\frac{1}{d})$, a saddle-node.
\item  if $\nu>\frac{2}{\sqrt{d}}$, then there exist three points:
$(F=0, G=0)$, a stable focus ($\nu<2$) or a stable node, otherwise,  and $(F=-\frac{\nu\pm\sqrt{\nu^2-\frac{4}{d}}}{2}, G=\frac{1}{d})$, a saddle and an unstable node.
\end{itemize}

We see that there are no equilibria in the domain $G>\frac{1}{d}$, hence there are no bounded trajectories in this region. If the motion on the plane $(G, F)$ starts from the point for which $G(0,r(0))>\frac{1}{d}$, then the phase trajectory rests in the half-plane $G>\frac{ 1}{d}$ and $G(t,r(t)) \rightarrow +\infty$ and $F(t,r(t)) \rightarrow -\infty$  for $ t \rightarrow t_*<\infty$. 



Thus, we study the problem only in the half-plane $G<\frac{1}{d}$, $F \in \mathbb{R}$, at which lies the point corresponding to the zero equilibrium state.

\subsection{Boundedness of the solution}\label{Bofs}
In the half-plane $ G<\frac{1}{d}$, $F \in \mathbb{R}$  system \eqref{FG4} has one equilibrium $(0,0)$. It corresponds to the stationary state $\bv=\bE=0$ and it is stable for any values of the parameters $d$ and $\nu>0$. Namely, as a linear analysis show,
\begin{itemize}
\item if $0<\nu<2 $ it is a stable focus;
\item if $\nu=2$ it is a stable degenerate node;
\item if  $\nu>2$
it is  a stable  node.
\end{itemize}

\begin{lem}\label{L1} There exists  $\delta>0$ such that if  the initial data $(F_0(r_0), G_0(r_0))$ belong to the $\delta$--neighborhood of the origin, $r_0\in \bar{\mathbb R}_+$, then any solution to
 \eqref{FG4} tends to zero exponentially as  $t\rightarrow +\infty$.
\end{lem}
The {\it proof} follows from the fact that
 $(0,0)$ is asymptotically stable for all choices of parameters $\nu$, $d$. $\square$

\begin{lem}\label{L2} Let $\Phi(G,F)=\Phi(G_0,F_0)$ be a closed phase curve corresponding to the solution of system \eqref{FG4} for $\nu=0$ with initial data $(F_0, G_0)$.  Then the phase curve corresponding to the solution of system \eqref{FG4} for  $\nu>0$ with initial data $(F_0, G_0)$ lies strictly inside the curve $\Phi(G,F)=\Phi(G_0,F_0)$.
\end{lem}

\proof
Let us construct the phase curve of \eqref{FG4} at $\nu=0$, i.e. the solution of
\begin{equation*}
 \frac{ \mathrm d F }{\mathrm d G}=-\frac{F^2+G}{F(1-d G)}.
 \end{equation*}
It implies
\begin{equation*}
\frac{ \mathrm d Z }{\mathrm d G}=-\frac{2}{1-d G}Z - \frac{2G}{1-d G}, \quad Z(G)=F^2.
 \end{equation*}
The solution is
\begin{eqnarray}\label{Int_d}
&&\Phi(G,F)=\frac{(d-2) F^2-2G+1}{(d-2)(1-d G)^{\frac{2}{d}}}=\Phi(G_0,F_0)=C_d,\\&& C_d=\frac{(d -2) F^2_0-2G_0+1}{(d-2)(1-d G_0)^{\frac{2}{d}}},\nonumber
 \end{eqnarray}
for $d\neq 2$, and
\begin{eqnarray}\label{Int_2}
&&\Phi(G,F)=\frac{2F^2+\ln(1-2G)(1-2G)+1}{2(1-2G)}=\Phi(G_0,F_0)=C_2, \\ &&C_2= \frac{2F^2_0+\ln(1-2G_0)(1-2G_0)+1}{2(1-2G_0)},\nonumber
\end{eqnarray}
 for $d=2$. As it was shown in \cite{R22_Rad}
 the curves given as \eqref{Int_2} and \eqref{Int_d} are bounded, they contain  the origin and intersect the axis $F=0$ in two points: $(G_-,0)$, $G_-<0$, and $(G_+,0)$, $G_+>0$, see the pictures in \cite{R22_Rad}.

Let us consider $V(t)=\Phi(G,F)$ as a Lyapunov function in the half-plane $ G<\frac{1}{d}$, $F \in \mathbb{R}$.
The derivative of $V(t)$ due  to system \eqref{FG4} is
\begin{eqnarray}\label{Lyap}
\frac{\mathrm d V}{\mathrm d t}=-\frac{2\nu F^2 }{(1-d G)^{\frac{2}{d}}}\le 0.
\end{eqnarray}
  If we denote  $(G(t), F(t))$ and  $(\bar G(t), \bar F(t))$ the point on the phase curve of \eqref{FG4} for $\nu>0$ and $\nu=0$, respectively, then
the distance  $|(G(t), F(t))|<|(\bar G(t), \bar F(t))|$, $t>0$, since $\frac{\mathrm d V}{\mathrm d t}=0$ if and only if $F=0$.  However $F=0$ does not solve \eqref{FG4}. $\square$

\bigskip

\begin{lem}\label{L2.1}
System \eqref{FG4}  has no limit cycles in the half-plane $G<\frac{1}{d}$, $F \in \mathbb{R}$.
\end{lem}
\proof
We use the Lyapunov function from Lemma \ref{L2} to
prove the absence of a limit cycle by contradiction. Assume that a limit cycle (a closed trajectory $\Gamma$) exists. Then it contains a stable equilibrium $(0,0)$ inside. We denote as $d(Y_1,Y_2)$ the  distance between points $Y_1(t)=(G_1(t),F_1(t)),\,Y_2(t)=(G_2(t),F_2(t))$.
For some initial point $Y(t_0)=(G_*,F_*)$ on $\Gamma$ there exists a time $t_1>t_0$ such that $Y(t_1)=Y(t_0)$ and, accordingly, $d(Y(t_0),Y(t_1))=0$. The curve $\Gamma$ contains  $(0,0)$ inside, so there are two points on this trajectory for which $F=0$, they are  $(G_+,0)$ and $ (G_-,0),$ $0<G_+<\frac{1}{d}, G_-<0$. At these points $\frac{\mathrm d V (G_+,0)}{\mathrm d t} =\frac{\mathrm d V (G_-,0)}{\mathrm d t}=0$. At other points of $\Gamma$ we have $\frac{\mathrm d V}{\mathrm d t}<0$. Then the function $V(t)$ does not increase along
$\Gamma$. Moreover, $\frac{\mathrm d V}{\mathrm d t}=0$ only at two points on $\Gamma$, therefore $V(t_1)-V (t_0)<0$ for $t_1>t_0$
and $d(Y(t_0),Y(t_1))>0$. Thus, we obtain a contradiction. $\square$

\bigskip

Thus, Lemmas \ref{L1},  \ref{L2} and \ref{L2.1} imply the following property of the phase trajectories:
\begin{lem}\label{L3} Let $d\ge 2$. Then the phase curves
  of system \eqref{FG4} are bounded in the half-plane $ G<\frac{1}{d}$, $F \in \mathbb{R}$ and tends to zero as  $t\rightarrow +\infty$.
\end{lem}

\begin{remark}
Lemma \ref{L3} is not valid for $d= 1$. Indeed, in this case  system \eqref{FG4} coincides with  system \eqref{Dlam} for the derivatives. As was shown in \cite{RChD20},   for any $\nu>0$ there exists a point on the phase plane such that the phase curve starting from this point goes to infinity as $t\to t_*<\infty$.
\end{remark}

\subsection{Study of the behavior of derivatives}
This section closely follows \cite{R22_Rad}, but for the convenience of the reader we give a sketch of the reasonings.

Let us denote $D=\Div\bV,$ $ \lambda=\Div{\bf  E}$.
Equations \eqref{4} imply
\begin{equation*}
\frac{\partial D}{\partial t}+(\bV \cdot \nabla D)=-D^2+2(d-1) FD- d(d-1) F^2-\lambda -\nu D,$$ $$\frac{\partial \lambda}{\partial t}+(\bV \cdot \nabla \lambda)= D(1-\lambda).
\end{equation*}
Along the characteristics given as  $\dot{r}=Fr$  the functions $D, \lambda$ obey
\begin{equation}\label{Dlam}
\dot D=-D^2+2(d-1) FD- d(d-1) F^2-\lambda -\nu D, \quad
\dot \lambda= D(1-\lambda).
\end{equation}
We introduce new variables
\begin{equation}\label{uvdef}
 u=D-d F, \qquad v=\lambda-d G.
 \end{equation}
 Systems \eqref{Dlam} and
\eqref{FG4} imply
\begin{eqnarray}
\label{uv}
\dot u =-u^2-2uF-v -\nu u,\quad \dot v=-uv+(1-d G)u-d F v.
\end{eqnarray}

System \eqref{uv} can be linearized my means of the Radon lemma (e.g. \cite{F}, \cite{Radon}).

\begin{theorem}\label{TR} [The Radon lemma]
\label{T2} A matrix Riccati equation
\begin{equation}
\label{Ric}
 \dot W =M_{21}(t) +M_{22}(t)  W - W M_{11}(t) - W M_{12}(t) W,
\end{equation}
 {\rm (}$W=W(t)$ is a matrix $(n\times m)$, $M_{21}$ is a matrix $(n\times m)$, $M_{22}$ is a matrix  $(m\times m)$, $M_{11}$ is a matrix  $(n\times n)$, $M_{12} $ is a matrix $(m\times n)${\rm )} is equivalent to the homogeneous linear matrix equation
\begin{equation}
\label{Lin}
 \dot Y =M(t) Y, \quad M=\left(\begin{array}{cc}M_{11}
 & M_{12}\\ M_{21}
 & M_{22}
  \end{array}\right),
\end{equation}
 {\rm (}$Y=Y(t)$  is a matrix $(n\times (n+m))$, $M$ is a matrix $((n+m)\times (n+m))$ {\rm )} in the following sense.

Let on some interval ${\mathcal J} \in \mathbb R$ the
matrix-function $\,Y(t)=\left(\begin{array}{c}{Q}(t)\\ {P}(t)
  \end{array}\right)$ {\rm (}${Q}$  is a matrix $(n\times n)$, ${P}$  is a matrix $(n\times m)${\rm ) } be a solution of \eqref{Lin}
  with the initial data
  \begin{equation*}\label{LinID}
  Y(0)=\left(\begin{array}{c}I\\ W_0
  \end{array}\right)
  \end{equation*}
   {\rm (}$ I $ is the identity matrix $(n\times n)$, $W_0$ is a constant matrix $(n\times m)${\rm ) } and  $\det {Q}\ne 0$ on ${\mathcal J}$.
  Then
{\bf $ W(t)={P}(t) {Q}^{-1}(t)$} is the solution of \eqref{Ric} with
$W(0)=W_0$ on ${\mathcal J}$.
\end{theorem}

Let us  \eqref{uv} as \eqref{Ric} with
\begin{eqnarray*}\label{M}
W=\begin{pmatrix}
  u\\
  v
\end{pmatrix},\quad
M_{11}=\begin{pmatrix}
  0\\
\end{pmatrix},\quad
 M_{12}=\begin{pmatrix}
  1 & 0\\
\end{pmatrix},\\
M_{21}=\begin{pmatrix}
 0 \\ 0
\end{pmatrix},\quad
M_{22}=\begin{pmatrix}
  - 2\,F - \nu & -1\\
  1-d\,G  & -d\,F\\
\end{pmatrix}.\\\nonumber
\end{eqnarray*}
Then according Theorem \ref{TR} the solition of \eqref{uv} is
\begin{eqnarray*}\
   \quad W(t)=\frac{P(t)}{Q(t)},
 \end{eqnarray*}
 where $P(t)=(p_1(t),p_2(t))^T$ and  $Q(t)$ solves the linear system
\begin{eqnarray}\label{LinQ4}
\left(\begin{array}{c}  Q\\  P
  \end{array} \right)^{\cdot}=M \left(\begin{array}{c}  Q\\  P
  \end{array} \right), \quad  M=\begin{pmatrix}
   0 & 1 & 0 \\
  0&-2F- \nu  & -1\\
  0&1-dG & - dF\\
\end{pmatrix},
 \end{eqnarray}
  subject to the initial data
\begin{eqnarray*}\label{cdQ4}
\left(\begin{array}{c}
  Q\\  P
  \end{array} \right)(0)=
  \left(\begin{array}{c}
   1\\  W_0
  \end{array} \right), \quad W_0= \left(\begin{array}{c}  u_0 \\  v_0
  \end{array} \right)=\left(\begin{array}{c}  \Div\bV_0-d F_0(r_0) \\  \Div  {\bf  E}_0-d G_0(r_0)
  \end{array} \right).
  \end{eqnarray*}
Since the vector-function $P(t)$ and the function $Q(t)$ are components of the solution of a linear system of differential equations with continuous coefficients, these functions do not go to infinity for any finite value of $t$. Hence the functions $u, v$ go to infinity along the characteristic starting from a point $ r_0 \in \mathbb R$ if and only if there exists a finite $t_*>0$ such that $Q(t_*,r_0)= 0$. Since $u=D-d F, v=\lambda-d G$ (see \eqref{uvdef}) and the functions $G, F$ are bounded, $D$ and $\lambda$  are bounded if and only if the functions $u, v$ are bounded. Since $D=\Div\bV,$ $ \lambda=\Div{\bf  E}$, due to the radial symmetry \eqref{sol_form}, we conclude that
 all derivatives of the solution of  \eqref{4}  are bounded if and only if the functions $u, v$ are bounded. Thus, the  conditions for the boundedness of  derivatives coincide with the conditions under which the function $Q(t)$ does not vanish for any finite $t$.

System \eqref{LinQ4} implies

\begin{eqnarray*}
Q(t)=1+\int\limits_0^t p_1(\tau) d\tau, \quad \dot p_1 =-(2F+\nu)p_1-p_2,\quad \dot p_2=(1-d G)p_1-d Fp_2.
\end{eqnarray*}
It follows
\begin{eqnarray*}
 \ddot p_1+((d+2) F+\nu)\dot p_1+( 2 \dot F +(1-d G)+d F(2F+\nu))p_1=0,\end{eqnarray*}
and, taking into account $ \dot F=-F^2-G-\nu F$,
\begin{eqnarray*}
 \ddot p_1+((d+2) F+\nu)\dot p_1+(2(d-1)F^2-(2+d)G+(d -2)\nu F +1)p_1=0.\end{eqnarray*}
We change \begin{eqnarray*}p_1(t)=H(t)e^{-\frac{\nu}{2}t}e^{-\frac{d+2}{2} \int\limits_0^t F(\tau) d\tau}\end{eqnarray*}
and obtain
\begin{eqnarray}\label{H}
\ddot H +JH=0,\end{eqnarray} with
\begin{eqnarray}\label{J}
J=1-\frac{1}{4}\nu^2-\frac{(d-2)(d-4)}{4}F^2+(d-2) \nu  F-\frac{(d+2)}{2}G.
\end{eqnarray}
Thus,
\begin{eqnarray}\label{QQQ}
Q(t)=1+\int\limits_0^t p_1(\tau) d\tau=1+\int\limits_0^t H(\tau)e^{-\frac{\nu}{2}\tau}e^{-\frac{d+2}{2} \int\limits_0^\tau F(\xi) d\xi} d\tau.\end{eqnarray}

Thus, for the boundedness of derivatives, it is necessary to require that for all $t>0$ condition
\begin{eqnarray}\label{intH0}
\int\limits_0^t H(\tau)e^{-\frac{\nu}{2}\tau}e^{-\frac{d+2}{2} \int\limits_0^\tau F(\xi) d\xi} d\tau>-1 \end{eqnarray}
holds.


It is easy to check that
\begin{equation}\label{H0}
H(0)=H_0=u_0, \quad \dot H(0)=H_1=\left(\frac{d -2}{2} F_0-\frac{\nu}{2}\right)u_0 -v_0.
\end{equation}

\bigskip

\section{Proof of Theorem \ref{T1}}

First of all, we notice that condition \eqref{T1CD} follows from
\begin{equation}\label{GFuv0}
\sup\limits_{r\in \overline{ \mathbb{R}}_+} |(r G_0(r), r F_0(r), u_0 (r), v_0(r))|<\varepsilon,
\end{equation}
where $|(x_1,\dots, x_k)|=\sqrt{x_1^2+\cdots+x_k^2}$, $k\in \mathbb N$.

Since we are interested in small values of $\nu$, we restrict ourselves to the case of $\nu<2$.


1. Let us fix $r\in \overline{ \mathbb{R}}_+ $.
 The matrix of linearization of the system of four equations \eqref{FG4}, \eqref{uv}, consists of two blocks
$
\begin{pmatrix}
  -\nu& -1\\
  1& 0
\end{pmatrix},
$
its complex conjugate eigenvalues are
$\lambda_{1,2}=-\frac{\nu \pm i h_1}{2},$ where
 $h_1={\sqrt{4-\nu^2}}$.
therefore the equilibrium  the origin on the phase space $G, F, u, v$ is asymptotically stable. This implies that
for any sufficiently small $\varepsilon$ - neighborhood of the origin there exists  $\delta(\varepsilon)<\varepsilon $ such that if
$|(r G_0, r F_0, u_0, v_0)|<\delta,$
then
$|(G, F, u, v)|<\varepsilon.$
Moreover,
 \begin{equation}\label{GFexp1}|r G, r F,u,v|\le C_1 e^{-\frac{\nu}{2} t},\quad C_1=\rm const,
 \end{equation}
  the constant  $C_1$ depends on $ \nu$ and $\varepsilon$, \cite{CL}, Ch.XIII, Sec.1,
  $\varepsilon=\varepsilon(\nu, r)$.

Thus, in  condition \eqref{GFuv0} we take $\varepsilon(\nu)=\inf\limits_{r\in \overline{ \mathbb{R}}_+} \varepsilon (\nu, r).$
The asymptotics \eqref{T1ass} follows immediately.

\bigskip

\section{Proof of Theorem \ref{T21}}


There is an alternative method to proof  Theorem \ref{T1}. Namely,  we can show that  for an arbitrary small $\nu>0$ there exists $\varepsilon (\nu)> 0$ such that if \eqref{GFuv0} holds,
then
\begin{eqnarray}
\label{intH}
\Big|\int\limits_0^t H(\tau)e^{-\frac{\nu}{2}\tau}e^{-\frac{d+2}{2} \int\limits_0^\tau F(\xi) d\xi} d\tau\Big|<1
\end{eqnarray}
for all $t>0$. This condition evidently implies \eqref{intH0}, therefore the derivatives of the considered solution are bounded, and the solution of \eqref{4}, \eqref{CD1} keeps smoothness.

However,  detailed estimates of the functions under the sign of integral gives us a possibility to obtain more or less practical sufficient condition that guarantees smoothness of solutions in terms of initial data.




1. Let us denote $S(t)=e^{-\frac{d+2}{2} \int\limits_0^t F(\xi) d\xi}$. Then $\dot S =-\frac{d+2}{2} F S $, and, as follows from \eqref{FG4},
\begin{eqnarray}\label{SG}
S= \Big|\frac{1-d G}{1-d G_0}  \Big|^{\frac{d+2}{2 d}},
\end{eqnarray}
therefore, due to \eqref{Lyap},
\begin{eqnarray}\label{M1}
0<M_-<S<M_+,\quad M_+=\Big|\frac{1-d G_-}{1-d G_0}  \Big|^{\frac{d+2}{2 d}}, \quad M_-=\Big|\frac{1-d G_+}{1-d G_0}  \Big|^{\frac{d+2}{2 d}},
\end{eqnarray}
where $G_-<0$ is the point, where the curve \eqref{Int_d} (or \eqref{Int_2}) intersects the axis $F=0$, see Lemma \ref{L2}.

2. To estimate $H(t)$ we use the following result \cite{Bellman} (Th.2, Ch.2):
all solution of the equation
\begin{equation*}\label{z}
\ddot z+(1+\varphi(t))z=0, \quad \int\limits^\infty |\varphi(\tau)| d\tau <\infty
\end{equation*}
are  bounded. Moreover, $z^2+ \dot z ^2\le  (y^2+{\dot y ^2}) \, e^{2\int\limits_0^t  |\varphi(\tau)| d\tau}$, where $y$ is a solution of $\ddot y+y=0$ such that $z(0)=y(0)$, $\dot z(0)=\dot y(0)$.
It implies
\begin{equation}\label{zy}
z^2\le  (y^2+{\dot y^2})\, e^{2\int\limits_0^t |\varphi(\tau)| d\tau}= (y^2(0)+{\dot y^2(0)})\, e^{2\int\limits_0^t |\varphi(\tau)| d\tau}.
\end{equation}

In \eqref{H} we can change the time variable as $t_1=h_1 t$,
 to obtain $\ddot H+J_1 H=0$, where $J_1=1+\varphi_1(t_1)$,
 \begin{equation*}
 \varphi_1(t_1)=\frac{1}{h_1}\left(-\frac{(d-2)(d-4)}{4}F^2+(d-2) \nu  F-\frac{(d+2)}{2}G\right).
 \end{equation*}
 From \eqref{GFexp1} we can conclude that
 $|\varphi_1(t_1)|\le {\rm const}\,  e^{-\frac{\nu}{2 h_1}t_1}$.
Since $\int\limits^\infty |\varphi_1(\tau)| d\tau  <\infty$,
then $|H(t_1)|$ is bounded for all initial data $H(0)$, $\dot H(0)$.
Now we go back to the time variable $t$ and use the notation of Theorem 2 for $\phi$ and $J_+$. Taking into account \eqref{zy} we get
 \begin{equation}\label{H1}
|H(t)|\le \sqrt{H_0^2+\frac{4  H_1^2}{4-\nu^2}} \,\, e^{\int\limits_0^\infty  |\phi(\tau)| d\tau}
 \end{equation}
 and
 \begin{eqnarray}\label{H2}
|H(t)|\le \sqrt{H_0^2+\frac{4 \dot H_0^2}{4-\nu^2}} \,\, e^{\int\limits_0^t  |\phi(\tau)| d\tau}\le
 \sqrt{H_0^2+\frac{4  H_1^2}{4-\nu^2}} \,\, e^{\left(J_+-1+\frac{\nu^2}{4}\right) T},
 \end{eqnarray}
for every $T>0$.

3. Note that due to Lemma \ref{L2}, for all $t>0$ the points $(G, F)$ lie inside the bounded curves \eqref{Int_2} or \eqref{Int_d}, therefore the maximal (positive) $G_+$ and minimal (negative) $G_-$ values of $G$, as well as maximum of  $F^2$, denoted as $F_+^2$,  can be found from the analytic expression for these curves. Therefore for every $(G_0, F_0)$ we can find $J_+=\rm const$ such that
$ J\le J_+$, where $J$ is given as \eqref{J}.

Estimates  \eqref{intH}, \eqref{M1}, \eqref{H1} imply condition \eqref{suffcond1}, whereas \eqref{intH}, \eqref{M1}, \eqref{H2}  imply condition \eqref{suffcond2}, if we  substitute \eqref{H0} and use the notation of Theorem 2.

We can only notice that  we do not need to know the value of $F_+$, which appear in $J_+$.  Indeed, for $d=2$ the expression for $J_+$ does not contain  $F_+$,
 for $d>2$ the value of $F_+$ can be found via $G$.

 Let us prove the latter statement.  At the maximum point of   \eqref{Int_d} we have $\dfrac{dF}{dG}=0$, i.e. $F^2=-G$. Therefore, the value of $G$, at which the extremum is reached, can be found from the equation
 \begin{equation*}
  -G=\frac{2G-1}{d-2}+(1-d G)^{\frac{2}{d}} C_d, \quad C_d=\frac{(d -2) F^2_0-2G_0+1}{(d-2)(1-d G_0)^{\frac{2}{d}}},
\end{equation*}
which solution is $G=\frac{1}{d} \left(1- ((d-2) C_d)^{\frac{d-2}{d}}  \right) $.
Thus,
 \begin{equation}\label{Fp}
  F_+=\frac{1}{d} \left(((d-2) C_d)^{\frac{d-2}{d}} -1   \right).
\end{equation}

4. Now we prove \eqref{suffcond3}. Let us denote as $H_*<0$ the value of $H(t)$ at the point $t_*$ of a negative minimum.  Assume that we know the estimate $H_*\le H_+^*<0$, then
\begin{eqnarray*}
\int\limits_0^t H(\tau)e^{-\frac{\nu}{2}\tau}e^{-\frac{d+2}{2} \int\limits_0^\tau F(\xi) d\xi} d\tau< \frac{2 H_+^* M_+}{\nu}<-1
\end{eqnarray*}
is a sufficient condition for the   blow-up.

Let $ H_+(t)$ be the solution to the Cauchy problem
\begin{eqnarray*}
\ddot { H}_+ + J_+  H_+=0, \quad  H_+(0)= H(0), \quad \dot{ H}_+(0)=\dot H(0).
\end{eqnarray*}
Indeed,
it is easy to check that
\begin{eqnarray*}
 \dfrac{d}{dt}\left( H^2+\frac{\dot H^2}{J_+}    \right) = \frac{2 (J_+-J)}{J_+} H \dot H.
\end{eqnarray*}

Since $H(0)\le 0$, $\dot H(0)<0$, then for $t\in (0,t_*)$ we have
 $H \dot H\ge 0$ and $H^2+\frac{\dot H^2}{J_+} \ge H^2(0)+\frac{\dot H^2(0)}{J_+}$, and in the point of minimum
$H_*^2 \ge H^2(0)+\frac{\dot H^2(0)}{J_+}\equiv (H_+^*)^2$. Note that $H(t) $ obtains its minimum on the semi-period of $H_+$, i.e. $t_*\le\frac{\pi}{\sqrt{J_+}}$.
Now it rests to substitute \eqref{H0}.

Thus, Theorem \ref{T21} is proved. 


\section{Proof of Theorem \ref{T22}}

Now we assume $\nu>2$ and fix $r_0\in \overline{\mathbb R}_+$.

1. The eigenvalues of the matrix of linearization of \eqref{FG4} are now real and negative:
$\lambda_{1,2}=-\frac{\nu \pm  h_2}{2},$ where
 $h_2={\sqrt{\nu^2-4}}$. Therefore (\cite{CL}, Ch.XIII, Sec.1)
 \begin{equation}\label{GFexp2}|(G,F)|\le C_2 e^{-\frac{\nu -  h_2}{2} t},\quad C_2={\rm const}>0.
 \end{equation}

2. We change the time as $t_2=\frac{h_2}{2} t$, and rewrite
 \eqref{H} as $\ddot H-J_2 H=0$, where $J_2=1+\varphi_2(t)$,
 \begin{equation*}
 \varphi_2(t)=-\frac{4}{h_2^2}\left(-\frac{(d-2)(d-4)}{4}F^2+(d-2) \nu  F-\frac{(d+2)}{2}G\right).
 \end{equation*}
The equation
  \begin{equation*}\label{u1}
u''-(1+\varphi_2(\tau))u=0, \quad  \int\limits^\infty |\varphi_2(\tau)| d\tau <\infty,
\end{equation*}
has two solution such that $u(t) \sim e^\tau$ and $u(t) \sim e^{-\tau}$ as $\tau \to \infty$ \cite{Bellman}.
Moreover, for $|\varphi_2|<1$ the solution  is non-oscillating and the equation has at most one  root for $t_2>0$.
Thus, \eqref{H} has two  non-oscillating solutions $H(t) \sim e^{\frac{h_2}{2} t}$ and $H(t) \sim e^{-\frac{h_2}{2} t}$ as $t \to \infty$.

3. Due to \eqref{intH0}, \eqref{M1}, it is enough to prove that
 \begin{eqnarray} \label{T2cond}
\Big|\int\limits_0^t H(\tau)e^{-\frac{\nu}{2}\tau} d\tau \Big|\to 0, \quad \nu\to\infty.
 \end{eqnarray}
 To this aim we perform twice the integration by parts to obtain
  \begin{equation*}
 \int\limits_0^t H(\tau)e^{-\frac{\nu}{2}\tau} d\tau = \Psi(H(t), \nu)+\int\limits_0^t H(\tau)e^{-\frac{\nu}{2}\tau} R(\tau) d\tau,
 \end{equation*}
 where
 \begin{equation*}\label{Psi}
 \Psi(H(t), \nu)=
 \frac{\nu}{2}(H(0)-H(t) e^{-\frac{\nu}{2} t})+\dot H(0)-\dot H(t) e^{-\frac{\nu}{2} t},
 \end{equation*}
  \begin{equation*}\label{R}
 R(t)=\frac{(d-2)(d-4)}{4}F^2+(d-2) \nu  F-\frac{(d+2)}{2} G.
 \end{equation*}
 Taking into account \eqref{H0}, we get
 \begin{equation*}\label{Psi1}
 \Psi(H(t), \nu)= \frac{d -2}{2} F_0 u_0 - v_0
 -\frac{\nu}{2} H(t) e^{-\frac{\nu}{2} t}-\dot H(t) e^{-\frac{\nu}{2} t}.
 \end{equation*}

4. Let us denote as $\bar H$ the solution of \eqref{H}, \eqref{H0} with $R=0$, which formally corresponds to $F=G=0$.
 This solution can be found explicitly as
  \begin{equation*}\label{Hbar}
  \bar H (t)=H(0)\cosh \frac{h_2 t}{2} +\frac{2 \dot H(0)}{h_2} \sinh \frac{h_2 t}{2},
 \end{equation*}
and
\begin{eqnarray*}
 &&\Psi(\bar H(t), \nu)=\frac{u_0}{h_2} e^{-\frac{\nu}{2} t} \sinh  \frac{h_2 t}{2}  \,  + \\\,&& \left(\frac{d -2}{2} F_0 u_0 - v_0\right)
 \left(1- e^{-\frac{\nu}{2} t}\left(\cosh \frac{h_2 t}{2}+ \frac{\nu}{h_2} \sinh \frac{h_2 t}{2}\right)\right).
 \end{eqnarray*}
It can be readily checked that for any fixed $t>0 $ as $\nu\to \infty$ we have
\begin{eqnarray*}
&&\frac{1}{h_2} e^{-\frac{\nu}{2} t} \sinh  \frac12 h_2 t= \frac{1}{\nu} + O \left(\frac{1}{\nu^2}\right), \\&& 1- e^{-\frac{\nu}{2} t}\left(\cosh \frac{h_2 t}{2}+ \frac{\nu}{h_2} \sinh \frac{h_2 t}{2}\right)= \frac{t}{\nu} + O \left(\frac{1}{\nu^2}\right),
 \end{eqnarray*}
 therefore $ \Psi(\bar H(t), \nu) \to 0$ as $\nu\to \infty$.

 5. Further we are going to prove that
  \begin{eqnarray}\label{w}
 H(t)=\bar H(t)+ O\left(\frac{1}{\nu}\right), \quad \dot H(t)=\dot{\bar H}(t)+ O(1), \quad \nu\to \infty, \quad  t>0.
  \end{eqnarray}

 Indeed, $w(t)=H(t)-\bar H(t) $ is the solution to the non-homogeneous problem
 \begin{eqnarray*}
\ddot w -\frac{h_2}{4} w= -R H, \qquad w(0)=\dot w(0),
 \end{eqnarray*}
 therefore, taking into account \eqref{GFexp2}, we have
  \begin{eqnarray*}
w(t)&=& \frac{1}{ h_2}\left(e^{-\frac{h_2 t}{2}} \int\limits_0^t R(\tau) H(\tau) e^{\frac{h_2 \tau}{2}} d\tau  - e^{\frac{h_2 t}{2}} \int\limits_0^t R(\tau) H(\tau) e^{-\frac{h_2 \tau}{2}} d\tau  \right)=\\ &&\frac{1}{h_2}\int\limits_0^t R(\tau) H(\tau) \sinh {{\frac{h_2 (\tau-t)}{2}}} d\tau = O\left(\frac{1}{\nu}\right), \quad \dot w(t)=  O(1), \quad \nu\to \infty.
 \end{eqnarray*}

 6. Now we show that for any fixed $t>0 $ as $\nu\to \infty$
 \begin{eqnarray}\label{RH}
 \int\limits_0^t H(\tau)e^{-\frac{\nu}{2}\tau} R(\tau)  d\tau=o\left(\int\limits_0^t H(\tau)e^{-\frac{\nu}{2}\tau} d\tau\right).
  \end{eqnarray}

   Indeed, \eqref{GFexp2} implies that there exists a constant $R_0>0$ such that $|R(t)|\le R_0 e^{-\frac{\nu -h_2}{2}t }$. Therefore
  \begin{eqnarray*}
&&\Big|\frac{1}{R_0} \int\limits_0^t H(\tau)e^{-\frac{\nu}{2}\tau} R(\tau)  d\tau - \int\limits_0^t H(\tau)e^{-\frac{\nu}{2}\tau} d\tau\Big|\le
\int\limits_0^t \Big|H(\tau)\Big| \Big|\frac{R(\tau)}{R_0}-1\Big|\,e^{-\frac{\nu}{2}\tau} d\tau \le \\
&&\int\limits_0^t |\bar H(\tau) + w(\tau)| |e^{-\frac{\nu -h_2}{2}\tau} -1|\,e^{-\frac{\nu}{2}\tau} d\tau = \int\limits_0^t \Big|\bar H(\tau) + O\left(\frac{1}{\nu}\right)\Big|\, e^{-\frac{\nu}{2}\tau}  d\tau\cdot O\left(\frac{1}{\nu}\right)=\\&&o\left(\frac{1}{\nu}\right) \to 0, \quad \nu\to\infty,
  \end{eqnarray*}
   what implies \eqref{RH}.

 7. Thus, for a fixed $t>0$ we have due to \eqref{w}
 \begin{eqnarray*}
 &&\Psi(H(t), \nu)=
 \frac{\nu}{2}(H(0)-(\bar H(t)+w(t)) e^{-\frac{\nu}{2} t})+\dot H(0)-(\dot {\bar H}(t)+\dot w(t)) e^{-\frac{\nu}{2} t}=\\
&& \Psi(\bar H(t), \nu) -(\dot w(t) + \frac{\nu}{2} w(t)) e^{-\frac{\nu}{2} t} \to 0,\quad \nu\to \infty.
 \end{eqnarray*}
Together with \eqref{RH}
 it implies \eqref{T2cond}.

  The asymptotic property \eqref{T2ass} can be proved as in Theorem 1.

  \bigskip

\section{Examples and numerical experiment}\label{Example}

The problem of blow-up or non-blow-up of solutions to  \eqref{4}, \eqref{CD1} for specific initial data and a specific coefficient $\nu$  can be solved numerically. Indeed, we solve system \eqref{FG4}, \eqref{H} for each $r$ and check the condition \eqref{intH0}.

For a numerical experiment, we select the initial data corresponding to a standard laser pulse \cite{CH18} (in fact, these are the only data that can be reproduced in the experiment).

Namely, we take
\begin{equation*}\label{DE}
\bE_0= a \, \bx \, \exp\left(-\frac{|\bx|^2}{2}\right), \qquad \bV_0=0, \qquad 0<a<\frac{1}{d}.
\end{equation*}
Since in the case $d>1$, $d\ne 4$, the solution blows up for all $a>0$ \cite{R22_Rad}, we are interested in studying the value of $\nu$ sufficient to guarantee a globally smooth solution.

For convenience, we write here the system, which follows from \eqref{FG4},  \eqref{H}, \eqref{QQQ},
\begin{eqnarray*}
  &&\dot F=-F^2-G-\nu F,\qquad \dot G=F-d FG, \qquad \dot {\mathcal F}=F, \\
  &&\dot H=Z,\qquad \dot Z=-J H,\qquad \dot Q=H \exp\left( -\frac{\nu}{2} \, t  - \frac{d+2}{2} \,\mathcal F \right),
\end{eqnarray*}
where $J$ is given as \eqref{J}, together with the initial data
\begin{eqnarray*}
  &&F(0)=0, \, G(0)={a}\, e^{-\frac{r^2}{2}}, \,  {\mathcal F}(0)=0,\, H(0)=0, \quad Z(0)={a\, r^2}e^{-\frac{r^2}{2}}, \, Q(0)=1,
  \end{eqnarray*}
  with
  $r\in \overline{\mathbb R}_+.$ Calculations are made by the fourth-fifth order Runge--Kutta--Felberg method (RKF45).
  \bigskip

 We present three types of numerical experiments. The behavior of $Q(t)$ essentially depends on the starting point of the characteristic $r$, so we have to do the calculations for the entire range of $r$ and choose the "worst" situation when $Q$ loses its positiveness. In fact, since the data decays exponentially at $r\to \infty$, it suffices to check the range of $r\in (0, 1.5)$.
\bigskip

  1. For $d=2$ we take $a=0.499$, that is, we fix $a$ very close to the possible upper value $\frac{1}{d}=0.5$ (recall that ${\rm div} \bE_0<1$). Computations show that $Q(t)$ attains its minimum over $r\in (0, 1.5)$ at $r\approx 0.03$, and the minimal value of $\nu$, sufficient to prevent a blow-up is approximately $0.9315$. Thus, this value is enough to prevent the breaking of (practically) all oscillations generated by a standard laser pulse. The graph of $Q(t)$ for this situation is presented in Fig.1, left.

  \bigskip

  2.  In \cite{RChD20} the physical value of $\nu$, based on the kinetic theory of gases, is computed, it is about $0.005$. If the temperature reasonably increases to keep the features of the model, $\nu$ increases up to $0.018$. In any case, this value is sufficiently less than the one, computed before. Let us study, what intensity of the laser pulse $a$ is allowed to generate oscillations that keep smoothness in the physical situation. If we take $\nu=0.005$, we can compute that the minimal value of $Q(t)$ attains approximately at $r=0.55$ and the maximal allowed value of $a$ is about $0.35$. For  $\nu=0.018$ the maximal value of $a$ is about $0.41$. The behavior of $Q(t)$ in the latter case is presented in Fig.1, right.

  \begin{figure}[htb]
\begin{minipage}{0.4\columnwidth}
\includegraphics[scale=0.25]{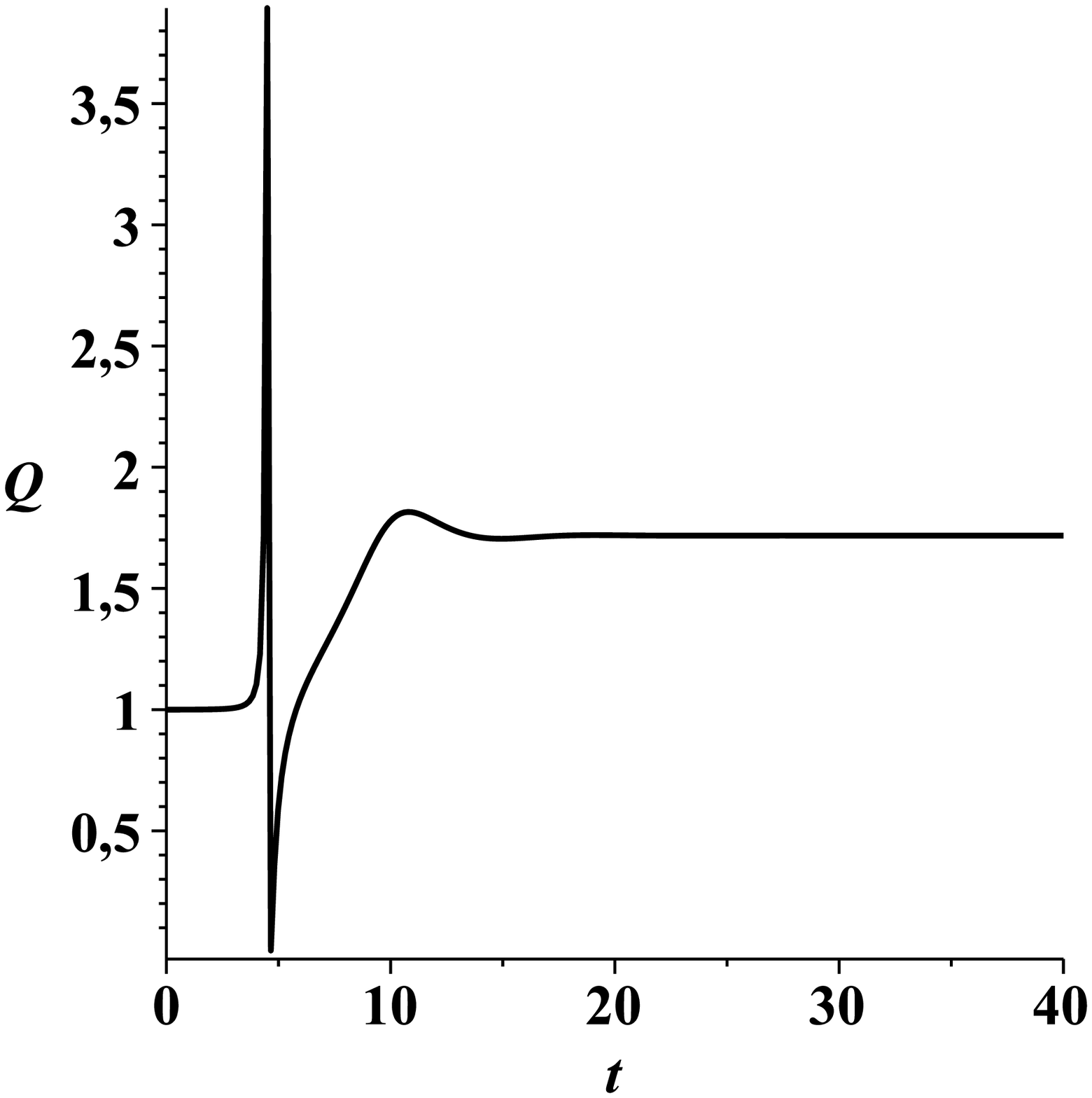}
\end{minipage}
\hspace{1.5cm}
\begin{minipage}{0.4\columnwidth}
\includegraphics[scale=0.25]{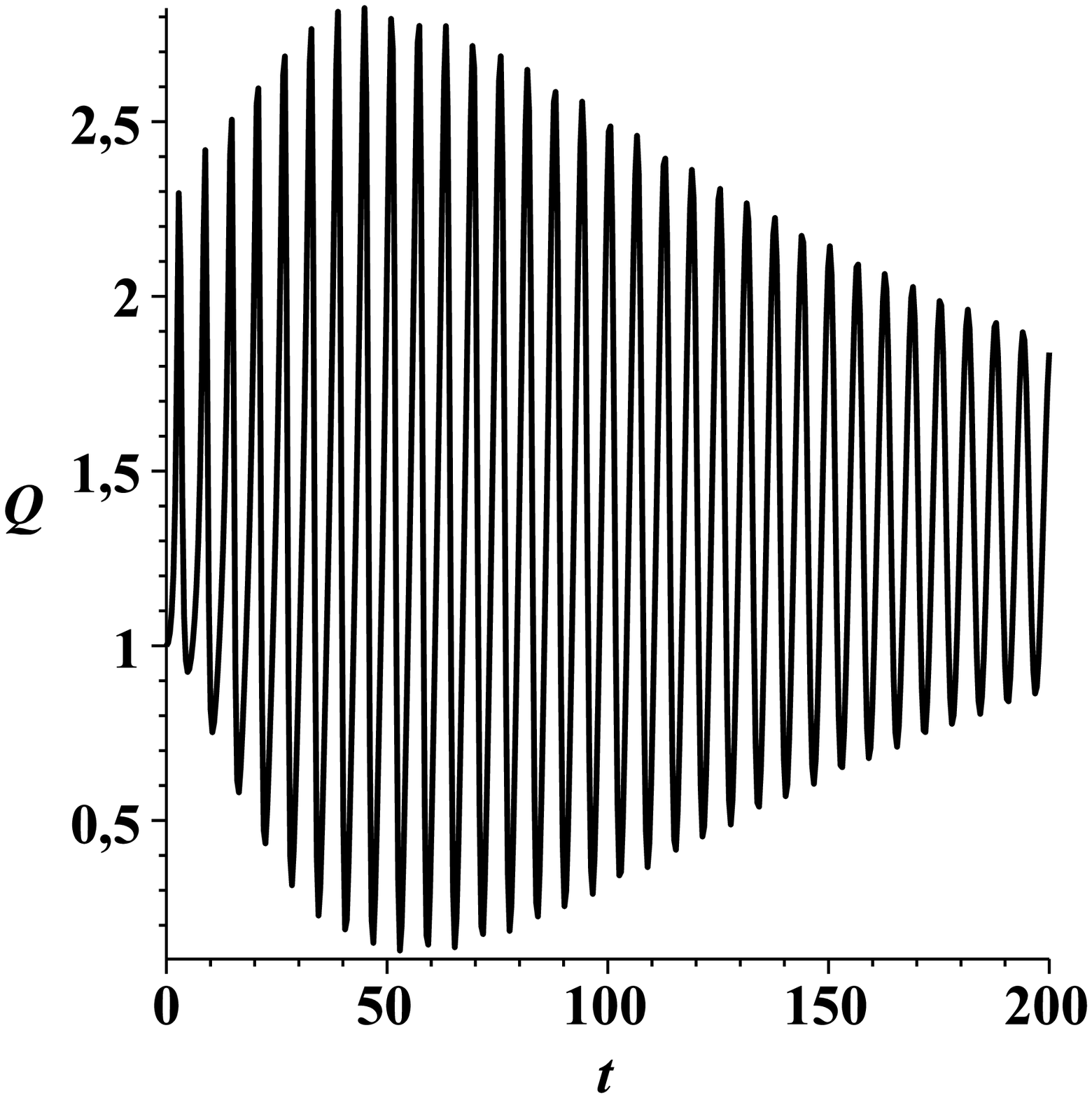}
\end{minipage}
\vspace{1cm}

\caption{The behavior of $Q(t)$ for $d=2$. Left: $a=0.499$, $\nu=0.9315$, $r=0.03$. Right: $a=0.41$, $\nu=0.018$, $r=0.35$ }\label{Pic4}
\end{figure}

  \bigskip

  3. It is also interesting that the properties of solution improve with the increase of dimension of space, $d$. Indeed, for $d=1$, to prevent breaking of oscillations generated by any standard laser pulse, one should take $\nu\ge 2$, see \cite{RChD20}, and for $d=2$, the value of $\nu$ is less, $0.9315$. Similar calculations carried out for $d=3$ show that $\nu=0.045$ is sufficient to prevent the blow-up of any ($a=0.33$) oscillations induced by the laser pulse.  For $\nu=0.005$, the maximum value of $a$ allowed for a smooth solution is $0.293$. Recall that for $d=3$ the intensity $a$ must be less than $\frac{1}{3}$.

\section{Discussion}

We proved that for axisymmetric multidimensional oscillations of  cold plasma the linear dumping, which corresponds to a constant coefficient of the frequency of collisions between particles $\nu$, serves as a mollifier. Moreover, Theorem 3 tells us that for an arbitrary initial pulse we can choose such a large coefficient $\nu$ that the solution  remains smooth for all $t>0$ and decay to the zero equilibrium state. However, this scenario does not make physical sense, since we cannot control the collision rate, which is relatively small ($\nu\ll 1$) according to the measurements.

The theoretical result of Theorem 1 is predictable.
Physicists know that small axisymmetric smooth deviations of the rest state persist in collisional media, see \cite{GFCA}, \cite{FCh19} for the cylindrical case and references therein.  They would be interested in the more or less exact size of the neighborhood of the zero equilibrium state corresponding to smooth solutions.
The criterion of smoothness in the terms of initial data can be obtained analytically for $d=1$, see \cite{RChD20}.
Theorem 2 gives some information about the lifetime of a smooth solution  for a fixed $\nu$.
However, this is not a criterion, but only sufficient conditions. The  condition \eqref{suffcond1} is more precise, but it is difficult to use in practice, since we do not know the analytical solution of \eqref{FG4}. Condition \eqref{suffcond2} is more rough than \eqref{suffcond1}, but more convenient, since we can  check arbitrary initial data \eqref{CD1} and decide what  lifetime  we can guaranty for the solution of the Cauchy problem    \eqref{4}, \eqref{CD1}.


Further, it should be noted that the constant collision frequency is only an assumption that simplifies the asymptotic analysis.
Actually $\nu$ is a function of density $n$.
It is shown in \cite{R_PhD} that in the case $d=1$ for $\nu = \nu_0 n^\gamma$, $\gamma>1$ each solution of the Cauchy problem is smooth for all initial data. A similar problem for the multidimensional case is completely open. It would be natural to expect that the form of $\nu(n)$ depends on $d$.

Another important problem is to study how  collisions between particles  affect solutions without radial symmetry.
The first approach to this difficult problem would be to study affine solutions for which
$(\bv, \bE)= (\mathbb F (t) \bx, \mathbb G(t) \bx)$, where $\mathbb F (t), \mathbb G (t)$ are matrices $(d \times d)$.
As shown in \cite{R22_Rad}, under the assumption of radial symmetry, such solutions are globally smooth. Nevertheless, as was recently proved \cite{RT22}, an arbitrarily small deviation from radial symmetry in the class of affine solution blows up, although the oscillation breaking mechanism is very subtle. The linearization  shows that the constant damping prevents the blow-up of asymmetric affine solutions. However, it is interesting to investigate whether this property holds for arbitrary asymmetric oscillations.

\section*{Acknowledgments}
O.R. was supported by the Russian Science Foundation (project No.23-11-00056) through RUDN University. M.D. was supported by the Moscow Center for
Fundamental and Applied Mathematics. This work does not have any conflicts of interest.


\end{document}